\documentclass[11pt,oneside]{amsart}
\usepackage{amsfonts}
\usepackage{amsmath}
\usepackage{amssymb}
\usepackage{mathrsfs}
\usepackage{geometry}
\usepackage{graphicx}
\usepackage{subfigure}
\usepackage{cite}
\usepackage{hyperref}
\usepackage{microtype}
\usepackage{enumerate}
\usepackage{enumitem}
\usepackage{tikz}
\usepackage{tikz-cd}
\usepackage{color}
\usepackage{ragged2e}
 \usepackage{float}
\allowdisplaybreaks[1]

\numberwithin{equation}{section}

\newcommand*{\circled}[1]{\lower.7ex\hbox{\tikz\draw (0pt, 0pt)%
    circle (.5em) node {\makebox[1em][c]{\small #1}};}}

\newcommand{\la}{\lambda}

\newcommand{\ga}{\gamma}

\newcommand{\ve}{\varepsilon}

\newcommand{\R}{\mathbb{R}}

\newcommand{\T}{\mathbb{T}}

\newcommand{\f}{\forall}

\newcommand{\n}[1]{\Vert #1\Vert }

\newcommand{\bbn}[1]{\Big\Vert #1 \Big \Vert }
\newcommand{\lr}[1]{\left\{ #1\right\} }
\newcommand{\lrc}[1]{\left[ #1\right] }
\newcommand{\lrs}[1]{\left( #1\right) }

\newcommand{\lra}[1]{\langle #1\rangle}

\newcommand{\wq}{\infty}
\newcommand{\pa}{\partial}

\newcommand{\ol}{\overline}

\theoremstyle{plain}
\newtheorem{theorem}{Theorem}[section]
\newtheorem{lemma}[theorem]{Lemma}
\newtheorem{proposition}[theorem]{Proposition}

\theoremstyle{definition}
\newtheorem{definition}[theorem]{Definition}

\theoremstyle{remark}
\newtheorem{remark}[theorem]{Remark}

\title[]{Uniqueness of global weak solutions for the
defocusing energy-critical NLS}

\author[X. Chen]{Xuwen Chen}
\address{Department of Mathematics, University of Rochester, Rochester, NY 14627, USA}
\email{xuwenmath@gmail.com}

\author[J. Holmer]{Justin Holmer}
\address{Department of Mathematics, Brown University, Providence, RI 02912, USA}
\email{justin\_holmer@brown.edu}

\author[S. Shen]{Shunlin Shen}
\address{School of Mathematical Sciences, University of Science and Technology of China, Hefei, 230026, China}
\email{slshen@ustc.edu.cn}

\begin{document}
\begin{abstract}
We study the defocusing energy-critical nonlinear Schr\"{o}dinger equation and establish global uniqueness for weak solutions.  The key idea comes from the weak uniqueness in the uniform-in-time high-frequency threshold (UTHT) class introduced in \cite{CH19,CH22unc}. We derive a new energy defect threshold (EDT) criterion from the UTHT criterion. Along the strong solution flow, energy conservation propagates the EDT condition for weak solutions and, combined with the global well-posedness theory for strong solutions, yields global weak uniqueness. As an application, we prove unconditional uniqueness for the quintic defocusing energy-critical Gross-Pitaevskii hierarchy.
\end{abstract}

\maketitle
\tableofcontents
\section{Introduction and main results}\label{sec:main}
We study the nonlinear Schr\"odinger equation (NLS)
\begin{align}\label{equ:NLS}
 \begin{cases}
 i\pa_{t} u=-\Delta u+F(u) ,\quad F(u)=\pm |u|^{p-1}u,\\
 u(t_{0},x)=u_{0}(x),
 \end{cases}
\end{align}
where $x\in \Lambda=\R^{d}$ or $\T^{d}$ and the signs $+$ and $-$ correspond to the defocusing and focusing cases, respectively.
In Euclidean spaces, the NLS \eqref{equ:NLS} is invariant under the scaling
\begin{align}\label{equ:NLS, scaling}
u_{\la}(t,x)=\la^{\frac{2}{p-1}}u(\la^{2}t,\la x),\quad \la>0,
\end{align}
which preserves the homogeneous Sobolev norm $\n{u_{0}}_{\dot{H}^{s_{c}}}$, where the critical scaling exponent is given by
\begin{align}\label{equ:scaling exponent}
s_{c}:=\frac{d}{2}-\frac{2}{p-1}.
\end{align}
In this paper, we focus on the cubic and quintic cases in the energy-critical setting $s_{c}=1$.

The global theory for the defocusing energy-critical NLS was developed through Bourgain's induction-on-energy argument for radial data \cite{Bou99} and the interaction Morawetz method of \cite{CKSTT08}, which removed the radial symmetry assumption. The concentration-compactness and rigidity method of Kenig--Merle \cite{KM06}
provided a further general approach to critical
evolution. We refer to \cite{Gri00,Tao05,RV07} for
related Euclidean developments and to
\cite{Bou93,HTT11,BD15,IP12,Kwak26,Yue21} for the
critical estimates and global theory on tori.
These works yield global strong flows in the
three-dimensional and four-dimensional cases considered here.

The purpose of this paper is to establish global
uniqueness for weak solutions of
\eqref{equ:NLS}. Given an initial datum $u_{0}\in H^{1}(\Lambda)$,  we call
 $u\in L^{\infty}(0,T;H^{1}(\Lambda))$ a weak solution if it satisfies the NLS \eqref{equ:NLS} in the sense of distributions.
For the energy-critical case, by Sobolev embedding,
\begin{align}
i\pa_{t}u\in L^{\infty}(0,T;H^{-1}),
\end{align}
and the standard weak-continuity criterion shows that every weak solution has a representative satisfying
\begin{align*}
 u\in W^{1,\infty}(0,T;H^{-1})\cap C([0,T];H^{1}_{\mathrm{weak}}).
\end{align*}
We always use this representative and interpret the initial condition in $H^{-1}$.

\begin{theorem}\label{theorem:weak uniqueness}
Let $\Lambda\in \lr{R^{3},\T^{3},\R^{4},\T^{4}}$, and let $I$ be a connected
open interval.
Let $u$ be a weak solution to the defocusing energy-critical NLS with
\begin{align*}
 u\in L^{\infty}_{\mathrm{loc}}(I;H^{1}(\Lambda)).
\end{align*}
Then, for every
$t_{0}\in I$,
\begin{align}\label{equ:identification}
 u(t)=\Phi(t-t_{0})u(t_{0}),\quad t\in I,
\end{align}
where $\Phi(t)$ denotes the global strong solution flow.
In particular, $u\in C(I;H^{1}(\Lambda))$ and $u$ enjoys the same space-time regularity as the global strong solution.
\end{theorem}

\vspace{1em}

\noindent\textbf{Quintic energy-critical Gross--Pitaevskii hierarchy}.
A direct application of Theorem~\ref{theorem:weak uniqueness} is unconditional uniqueness for the Gross--Pitaevskii
(GP) hierarchy. These hierarchies arise in the derivation of NLS from quantum many-body dynamics.
A basic problem is to justify the effective nonlinear evolution of a Bose--Einstein condensate,
starting from the linear Schr\"odinger equation for the many-particle wave function. In the hierarchy approach,
one studies reduced density matrices through the Bogoliubov--Born--Green--Kirkwood--Yvon
(BBGKY) hierarchy and identifies their limits
with solutions of a GP hierarchy.

Effective three-body interactions lead to a
quintic nonlinearity, which is energy-critical
in three dimensions. The corresponding
defocusing GP hierarchy is
\begin{equation}\label{equ:gp hierarchy,quintic}
 i\pa_{t}\ga^{(k)}=\sum_{j=1}^{k}[-\Delta_{x_{j}},\ga^{(k)}]
 +b_{0}\sum_{j=1}^{k} B_{j;k+1,k+2}\gamma^{(k+2)},
 \quad k\geq1, \quad b_{0}>0,
\end{equation}
where $b_{0}>0$ is a coupling constant.
For simplicity, we set $b_{0}=1$. The collision operator
$B_{j;k+1,k+2}$ is defined by
\begin{align}\label{equ:kernel,gp}
B_{j;k+1,k+2}\gamma^{(k+2)}:=
\operatorname{Tr}_{k+1,k+2}[\delta(x_{j}-x_{k+1})\delta(x_{j}-x_{k+2}),\ga^{(k+2)}].
\end{align}
 If $u$ solves the quintic NLS, then
\begin{align*}
\gamma^{(k)}(t)=|u(t)\rangle\langle u(t)|^{\otimes k}, \quad k\geq1,
\end{align*}
solves the quintic GP hierarchy \eqref{equ:gp hierarchy,quintic}. In other words, the GP hierarchy contains the NLS flow through its factorized solutions. Moreover, the weak solutions of the NLS are connected with the weak solutions of the GP hierarchy.

We denote by $\mathcal{L}_{k}^{1}$ the space of trace class operators on the
$k$-particle symmetric Hilbert space $L_{\text{sym}}^{2}(\Lambda^{k})$,
and use the notation
\begin{align*}
S^{(\pm1,k)}=\prod_{j=1}^k\langle\nabla_{x_j}\rangle^{\pm1},\quad \lra{\nabla}=(1-\Delta)^{\frac{1}{2}}.
\end{align*}
We then give the definition of weak solutions of the quintic GP hierarchy.
\begin{definition}\label{def:weak solution,gp}
A sequence $\Gamma=\{\gamma^{(k)}\}_{k=1}^{\infty}$ with
\begin{align*}
 \gamma^{(k)}\in L^{\infty}(0,T;\mathcal L_{k}^{1}),  \quad k\geq1,
\end{align*}
is called a weak solution of \eqref{equ:gp hierarchy,quintic} on $[0,T]$ if
$\gamma^{(k)}(t)$ is symmetric and nonnegative for every
$t$ and the following conditions hold.
\begin{enumerate}[wide=0pt, label=(\arabic*), align=left]
\item
For every $t\in[0,T]$, there is a $U(1)$-invariant
Borel probability measure $\mu_t$ on the closed unit ball of
$L^2(\Lambda)$ such that
\begin{equation}\label{equ:gamma,k,t,mut}
 \gamma^{(k)}(t)=\int(|\phi\rangle\langle\phi|)^{\otimes k}
d\mu_t(\phi),\quad k\geq1.
\end{equation}

\item
There is a constant $C_0>0$, independent of $k$, such that
\begin{equation}\label{equ:kinetic,bound}
 \operatorname*{sup}_{0\leq t\leq T}
 \bigl\|S^{(1,k)}\gamma^{(k)}(t)S^{(1,k)}
 \bigr\|_{\mathcal L_k^1}\le C_0^{2k},
 \quad k\geq1.
\end{equation}

\item For every $t\in[0,T]$ and $k\geq1$, the hierarchy satisfies
\begin{align*}
\gamma^{(k)}(t)
=&\gamma^{(k)}(0)
-i\int_{0}^{t}
\left(
\sum_{j=1}^k[-\Delta_{x_j},\gamma^{(k)}(\tau)]
+b_0\sum_{j=1}^k
B_{j;k+1,k+2}\gamma^{(k+2)}(\tau)
\right)\,d\tau,
\end{align*}
where, after conjugation by $S^{(-1,k)}$, the identity holds
in $\mathcal{L}_{k}^{1}$.

\item For every $k\geq1$, the map
\begin{align*}
t\longmapsto S^{(-1,k)}\gamma^{(k)}(t)S^{(-1,k)}
\end{align*}
is weak-* continuous in $\mathcal{L}_{k}^{1}$ on $[0,T]$, where the weak-* topology is induced by testing against compact operators.
\end{enumerate}
\end{definition}
\begin{remark}
The quantum de Finetti theorem implies condition (1) under either of the following hypotheses (see \cite[Definition~4.5 and Theorem~4.6]{CH19} for further details).
\begin{enumerate}[wide=0pt, label=(\arabic*), align=left]
\item[(a)]
The hierarchy is admissible, namely,
\begin{align*}
 \operatorname{Tr}\gamma^{(k)}(t)=1,
 \quad
 \operatorname{Tr}_{k+1}\gamma^{(k+1)}(t)=\gamma^{(k)}(t),
 \quad k\geq1.
\end{align*}
\item[(b)]
There are nonnegative bosonic $N$-particle density matrices $\gamma_N^{(N)}(t)$ with trace one such that, along a common
sequence $N\to\infty$,
\begin{align*}
\operatorname{Tr}_{k+1,\ldots,N}\gamma_N^{(N)}(t)\rightharpoonup^*\gamma^{(k)}(t)
 \quad\text{in $\mathcal{L}_{k}^{1}$ for every fixed $k$}.
\end{align*}
Here, weak-$*$ convergence is tested against compact operators.
\end{enumerate}
In case (a), $\mu_{t}$ is supported on the $L^{2}$ unit sphere.
In case (b), it is supported on the closed unit ball.
\end{remark}

\begin{theorem}[Uniqueness for the GP hierarchy]\label{thm:hierarchy}
Let $\Lambda\in\lr{\R^{3},\T^{3}}$.
For each initial hierarchy $\Gamma_{0}$, there is at most one weak solution
$\Gamma=\{\gamma^{(k)}\}_{k=1}^{\infty}$ of
\eqref{equ:gp hierarchy,quintic} on $[0,T]$ in the sense of
Definition \ref{def:weak solution,gp} satisfying $\Gamma(0)=\Gamma_{0}$.
\end{theorem}

\subsection{Key ideas and outline of the proof}
The key idea of the uniqueness argument for weak solutions comes from the uniqueness analysis of the GP hierarchy itself. Specifically, a crucial ingredient is the uniqueness established for the class of weak solutions introduced in \cite{CH19,CH22unc} under the uniform-in-time high-frequency threshold (UTHT) condition, which reads
\begin{align}
\n{u(t)}_{L_{t}^{\infty}([0,T_{0}];H^{1})} \leq C_{0},\quad
\n{P_{>M} u(t)}_{L_{t}^{\infty}([0,T_{0}];H^{1})} \leq \eta(C_{0})
\end{align}
for some $M\geq 1$,
where $\eta(C_{0})$ is a threshold that depends only on the $H^{1}$ bound $C_{0}$.

We briefly review significant advances along this line of research, which served as important inspiration for establishing the uniqueness of the class of weak solutions under the UTHT condition and Theorems \ref{theorem:weak uniqueness} and \ref{thm:hierarchy}.

One of the main purposes of proving uniqueness for the GP hierarchy is to establish the rigorous derivation of the NLS from quantum many-body dynamics. The first major breakthrough is due to Erd\H{o}s, Schlein, and Yau, who around 2005 rigorously derived the three-dimensional cubic defocusing NLS in a series of fundamental papers \cite{ESY06,ESY07,ESY09,ESY10}. A crucial ingredient was an
 $H^{1}$-type unconditional uniqueness theorem for the cubic GP hierarchy, proved in \cite{ESY07} through a sophisticated Feynman graph analysis. These groundbreaking papers have motivated much subsequent work.

Subsequently in 2007, Klainerman and Machedon \cite{KM08} proved a different uniqueness theorem under an additional a priori space-time norm assumption. Drawing on \cite{ESY07,KM93}, they introduced space-time collapsing estimates together with a different combinatorial scheme, now known as the Klainerman-Machedon (KM) board game argument. This argument has been employed in numerous works on the uniqueness of the GP hierarchy and the derivation of the NLS, including \cite{CP11,CP14der,CH16,CH16cor,CH22,CS14,GSS14,HS16,KSS11,She21,Soh15}. Moreover, although originally tailored to the NLS, it later provided essential insights into the Boltzmann hierarchy, as demonstrated for instance in \cite{AMPT25,CDP19,CH23der,CSZ26}.

In 2013, T.~Chen--Hainzl--Pavlovi\'c--Seiringer \cite{CHPS15} introduced the quantum de Finetti theorem into
hierarchy analysis, reducing estimates for density
matrices to multilinear estimates. This gave a simplified proof of
unconditional uniqueness. Their method, integrating the KM board game argument, the quantum de Finetti theorem, and multilinear Sobolev estimates, is robust in handling unconditional uniqueness problems of GP hierarchies.  For subsequent applications, see for instance \cite{HTX15,HTX16}.

Somewhat unexpectedly, the uniqueness analysis of the GP hierarchy, although seemingly more difficult than the corresponding analysis for the NLS, began to yield novel results for the NLS itself, starting with the work of Herr and Sohinger \cite{HS19} and of X. Chen and Holmer \cite{CH19}.
In \cite{HS19}, Herr and Sohinger generalised the usual
Sobolev multilinear estimates to Besov spaces and obtained new unconditional-uniqueness results
for the cubic NLS on $\T^{d}$, covering the full scaling-subcritical regime for $d\geq 4$. (See also \cite{Kis21} using NLS analysis in the scaling-subcritical regime.)

In the treatment of the
energy-critical GP hierarchy, Chen--Holmer \cite{CH19} introduced hierarchically uniform frequency localization (HUFL), and established a new $H^{1}$-type uniqueness theorem for the $\T^{3}$ quintic energy-critical GP hierarchy.
 For factorized hierarchies, imposing this localization uniformly
in time leads to the uniform-in-time
frequency localization (UTFL) property, which implies the
unconditional uniqueness result for the $\T^{3}$ quintic energy-critical NLS. In subsequent work \cite{CH22unc,CSZ22}, new combinatorics were worked out to enable the application of
$U$-$V$ multilinear estimates, which plays a key role in completely resolving the unconditional uniqueness problems for cubic and quintic energy-critical and energy-supercritical NLS on $\R^{d}$ and $\T^{d}$.

The decisive point for weak solutions is that the hierarchy method for uniqueness depends on only one high-frequency threshold. Thus, the full UTFL
condition can be replaced by a uniform bound on the high-frequency estimate at that threshold.
This is the UTHT criterion of Chen--Holmer
\cite[Theorem 5.12]{CH19} and \cite[Corollary~3.2]{CH22unc}, which applies beyond
the strongly continuous solution class $C([0,T];H^{1})$.

\vspace{1em}
\noindent \textbf{Outline of the proof}

In Section \ref{section:Conditional weak uniqueness and simplified proof}, we present a new and simplified proof of uniqueness in the UTHT weak solution class.
Although the argument is totally independent of the hierarchy method, its underlying ideas stem from our understanding of the hierarchy analysis. The first key observation is the local propagation of the UTHT property under the strong solution flow. We consider the class
\begin{align}
\mathcal{S}(C_{0},M,\eta):=\lr{f\in H^{1}\big| \ \n{f}_{H^{1}}\leq C_{0},\ \n{P_{>M}f}_{H^{1}}\leq \eta},
\end{align}
where $\eta$ is a threshold.
 We show that, for a local time $T(C_{0},M,\eta)>0$, the strong solution flow $\Phi(t)$ preserves the UTHT property according to
\begin{align}\label{equ:strong flow,UTHT,intro}
\Phi(t)\mathcal{S}(C_{0},M,\eta)\subset \mathcal{S}(2C_{0},M,2\eta),\quad t\in[0,T].
\end{align}
This observation in fact originates from the local propagation of the HUFL property under the strong solution flow introduced in \cite{CH19} in the analysis of the GP hierarchy.
Based on \eqref{equ:strong flow,UTHT,intro}, and via high-low frequency $U$-$V$ multilinear estimates, we then establish the uniform stability under the strong solution flow
\begin{align}\label{equ:stability,intro}
\sup_{t\in[0,T]}
\|\Phi(t)f-\Phi(t)g\|_{H^{-1}}
\leq2\|f-g\|_{H^{-1}},\quad \f f,g\in \mathcal{S}(C_{0},M,\eta).
\end{align}
With these stability estimates in hand, we conclude the weak uniqueness within the UTHT class via a simplified strong-weak uniqueness argument.

A crucial feature of the UTHT condition is that uniqueness requires only a single threshold. On the other hand, for the defocusing case, the high-frequency kinetic energy is closely related to the total energy. Therefore, starting from the uniqueness of the UTHT weak solution class, we provide a new energy defect threshold criterion for a class of weak solutions. More precisely,
in Section \ref{section:Local regularity below the energy defect threshold},
for a weak solution $u(t)$ with a weakly $H^{1}$ continuous representative, we define
\begin{align}\label{equ:defect,intro}
D(u,t_{0})=\lim_{\ve \to 0}\operatorname*{esssup}_{|t-t_{0}|<\ve}E(u(t))-E(u(t_{0})),
\end{align}
where the standard energy is given by
\begin{equation}\label{equ:intro,energy conservation}
 E(u)=\frac{1}{2}\n{\nabla u}_{2}^{2} +\frac{1}{p+1}\n{u}_{p+1}^{p+1}.
\end{equation}
We say that $u(t)$ satisfies the energy defect threshold (EDT) condition at $t_{0}$ if
$D(u,t_{0})<\delta$, where $\delta>0$ depends
on the a priori $H^{1}$ bounds and the UTHT threshold $\eta$.
Based on the convexity of the defocusing energy and the weak
$H^{1}$ continuity, we have
\begin{align*}
\limsup_{t\to t_{0}}\n{\nabla u(t)-\nabla u(t_{0})}_{L^{2}}^{2}\leq 2D(u,t_{0}),
\end{align*}
which suffices to ensure the UTHT property near $t_{0}$,
and hence yields local weak-strong uniqueness in the EDT class near $t_{0}$.

For a weak solution, the remaining problem is to obtain the EDT condition at some starting time and to retain it under continuation, yet energy conservation is not available at this stage.
In Section \ref{section:Strong-weak uniqueness}, we address this problem by choosing the starting time from the energy profile. Let $I_{0}\Subset I$ be a bounded
open interval and set
\begin{align*}
 A=\operatorname*{ess\,sup}_{t\in I_{0}}E(u(t)).
\end{align*}
We may take $t_{0}\in I_{0}$ satisfying
$E(u(t_{0}))>A-\delta$, so that the EDT condition is available at $t_{0}$.
Applying the local weak-strong uniqueness within the EDT class, we deduce that near $t_{0}$
the weak solution in fact coincides with the strong solution, namely
\begin{align}\label{equ:ut=vt,t0,intro}
u(t)=\Phi(t-t_{0})u(t_{0}),\quad \f t\in(t_{0}-\ve_{0},t_{0}+\ve_{0}).
\end{align}
We observe that the strong solution $\Phi(t-t_{0})u(t_{0})$ with initial datum $u(t_{0})$ conserves energy
\begin{align}\label{equ:energy conservation,vt,ut0,intro}
E(\Phi(t-t_{0})u(t_{0}))=E(u(t_{0})).
\end{align}
In view of \eqref{equ:ut=vt,t0,intro} and \eqref{equ:energy conservation,vt,ut0,intro}, the weak solution $u(t)$ satisfies the EDT property in a neighborhood of $t_{0}$. By a standard bootstrap argument, this identification extends throughout the interval on which the strong solution exists. Finally, relying on the well-known global well-posedness theory \cite{Bou99,CKSTT08,HTT11,IP12,RV07,Yue21}, we obtain global weak uniqueness.

The uniqueness theory developed above also reaches beyond the NLS itself. As a direct application of Theorem \ref{theorem:weak uniqueness}, in Section \ref{section:Uniqueness for the Energy-critical GP hierarchy} we prove unconditional uniqueness for the energy-critical quintic GP hierarchy. The bridge is provided by Ammari, Liard, and Rouffort \cite{ALR20}, who established a principle transferring uniqueness of weak NLS solutions to uniqueness for the associated GP hierarchy. Thanks to this uniqueness principle, the weak uniqueness theorem proved for the NLS extends to the hierarchy, thereby resolving the unconditional uniqueness problem at the hierarchy level.

\section{A UTHT criterion for weak uniqueness and a simplified proof}\label{section:Conditional weak uniqueness and simplified proof}
In this section, we begin by recalling the established uniqueness result for the UTHT weak solution class.
\begin{proposition}[UTHT criterion, {\cite{CH19,CH22unc}}]\label{lemma:weak uniqueness,threshold}
For every $C_{0}\geq 1$, there exists a threshold $\eta(C_{0})>0$ with the following property.
Given an initial datum $u_{0}$,
there is at most one weak solution $u(t)\in C([0,T_{0}];H_{\mathrm{weak}}^{1})$ to \eqref{equ:NLS} for the energy-critical cubic and quintic cases, satisfying the UTHT condition, namely the two requirements below.
\begin{enumerate}[wide=0pt, label=(\arabic*), align=left]
\item The solution $u(t)$ satisfies an $H^{1}$ energy bound
\begin{align*}
\n{u(t)}_{L_{t}^{\infty}([0,T_{0}];H^{1})} \leq C_{0} .
\end{align*}
\item There exists some frequency $M$ such that
\begin{align*}
\n{P_{>M} u(t)}_{L_{t}^{\infty}([0,T_{0}];H^{1})} \leq \eta(C_{0}).
\end{align*}
\end{enumerate}
Moreover, every such weak solution satisfies
\begin{align}\label{equ:identity,entire,prop}
u(t)=\Phi(t)u(0),\quad t\in [0,T_{0}].
\end{align}
\end{proposition}

Our goal is to provide a new and simplified proof of Proposition \ref{lemma:weak uniqueness,threshold}.
The underlying idea is local propagation of the UTHT property under the strong solution flow, which in fact originates from the local propagation of the HUFL property for the GP hierarchy in \cite{CH19}. Then, based on this propagation, we establish
the low-regularity uniform stability estimate under the strong solution flow.
The local propagation and uniform stability estimate are key to the simplified proof and may have broader applications, as explained in Remark \ref{remark:uniform stability}.

We first prove the proposition on $\mathbb{T}^3$ and $\mathbb{T}^4$ using uniform $H^{-1}$ stability of the strong flow. We then give a direct proof on $\mathbb{R}^3$ and $\mathbb{R}^4$ using endpoint Strichartz estimates.
\subsection{The periodic domains}
We recall the solution space $X^{s}([0,T])$,
 which is built on the $U$-$V$ spaces introduced by Koch-Tataru \cite{KT05,KT07}.
 We omit the definition and refer to \cite{HTT11,HTT14,IP12,KV16} for further details. The estimates used below are summarised as follows:
\begin{align}
\n{u}_{L_{t}^{\wq}H_{x}^{s}}\lesssim &\n{u}_{X^{s}},\label{equ:u-v and sobolev}\\
 \n{e^{it\Delta}f}_{X^{s}}\leq&\n{f}_{H^{s}}, \label{equ:u-v and sobolev, special case}\\
 \bbn{\int_{0}^{t}e^{i(t-\tau)\Delta}u(\tau,\cdot)d\tau}_{X^{s}([0,T])}\lesssim & \n{u}_{L_{t}^{1}H^{s}}.\label{equ:X,duhamel}
\end{align}

We then establish the uniform $H^{-1}$ stability estimate under the strong solution flow.
\begin{lemma}[Local propagation and uniform stability]\label{lemma:stability,strong flow,UTFLS}
Let
\begin{align}
C_{0}\geq1,\quad M\geq1,\quad 0<\eta\ll C_{0}^{-(p-2)},
\end{align}
and
\begin{align}
\mathcal{S}(C_{0},M,\eta):=\lr{f\in H^{1}\big| \ \n{f}_{H^{1}}\leq C_{0},\ \n{P_{>M}f}_{H^{1}}\leq \eta}.
\end{align}
There is a time $T=T(C_{0},M,\eta)\in(0,1)$ such that,
for $f\in \mathcal{S}(C_{0},M,\eta)$, the local strong solution $\Phi(t)f$ exists on $[0,T]$ and satisfies
\begin{align}\label{equ:local propagation,strong flow}
\|\Phi(\cdot)f\|_{X^{1}([0,T])}\leq2C_{0},
\quad
\|P_{>M}\Phi(\cdot)f\|_{X^{1}([0,T])}\leq2\eta.
\end{align}
Furthermore, if $f$ and $g$ both belong to $\mathcal{S}(C_{0},M,\eta)$, then
\begin{align}\label{equ:stability}
\sup_{t\in[0,T]}
\|\Phi(t)f-\Phi(t)g\|_{H^{-1}}
\leq2\|f-g\|_{H^{-1}}.
\end{align}

\end{lemma}

\begin{proof}
We first treat the cubic case on $\T^{4}$,
and record the estimates from
\cite[Lemmas 2.1 and 2.2]{CH22unc},
\begin{align}\label{equ:trilinear,high,low}
&\bbn{
\int_{0}^{t} e^{i(t-s)\Delta}
\widetilde{v}_{1}(s)\widetilde{v}_{2}(s)\widetilde{v}_{3}(s)ds}_{X^\sigma}\\
\lesssim&
\n{v_{1}}_{X^\sigma}
\lrs{T^{\frac{1}{7}}M^{\frac{3}{5}}\n{P_{\leq M}v_{2}}_{X^{1}}
+\n{P_{>M}v_{2}}_{X^{1}}}\n{v_{3}}_{X^{1}},\notag
\end{align}
where $\sigma\in\{-1,1\}$ and $\widetilde{v}_{j}$ denotes either $v_{j}$ or $\overline{v_{j}}$.

Let us define the nonlinear map
\begin{align*}
\mathcal{T}_{f}[u](t):= e^{it\Delta}f
-i\int_{0}^t e^{i(t-s)\Delta}F(u(s))ds,
\end{align*}
for fixed $f\in \mathcal{S}(C_{0},M,\eta)$ on the closed subset of $X^{1}([0,T])$ given by
\begin{align}
B(C_{0},\eta):=\lr{u\in X^{1}([0,T])\big|\
\n{u}_{X^{1}}\leq 2C_{0},\  \n{P_{>M}u}_{X^{1}}\leq 2\eta}.
\end{align}

We decompose the nonlinearity as
\begin{align*}
|u|^{2}u=&|P_{\leq M}u+P_{>M}u|^{2}(P_{\leq M}u+P_{> M}u)
=|P_{\leq M}u|^{2}P_{\leq M}u+F_{\text{high}}(u),
\end{align*}
where $F_{\text{high}}(u)$ denotes the sum of the terms involving at
least one factor of $P_{>M}u$.

For the low-frequency term, by Sobolev embedding and interpolation, we have
\begin{align*}
\n{|v|^2 v}_{H^1} \lesssim\n{v}_{L^{\infty}}^{2}\n{v}_{H^{1}}  \lesssim\n{v}_{H^1}\n{v}_{H^{3}}^{2} .
\end{align*}
Applying the above estimate, together with $X^{s}$ estimates \eqref{equ:u-v and sobolev}--\eqref{equ:X,duhamel} and Bernstein inequalities, we derive
\begin{align}\label{equ:low frequency,X1}
&\bbn{
\int_{0}^t e^{i(t-s)\Delta}\lrs{|P_{\leq M}u|^{2}P_{\leq M}u}(s)ds
}_{X^{1}}\\
\lesssim& T\n{P_{\leq M}u}_{X^{3}}^{2}\n{P_{\leq M}u}_{X^{1}}\lesssim TM^{4}C_{0}^{3}.\notag
\end{align}

For the term $F_{\text{high}}(u)$, we use estimate \eqref{equ:trilinear,high,low} to get
\begin{align}\label{equ:high frequency,X1}
&\bbn{ \int_{0}^{t} e^{i(t-s)\Delta}F_{\text{high}}(u(s))ds}_{X^{1}}\\
\lesssim& \lrs{T^{\frac{1}{7}}M^{\frac{3}{5}}\n{P_{\leq M}u}_{X^{1}}
+\n{P_{>M}u}_{X^{1}}}\n{u}_{X^{1}}\n{P_{>M}u}_{X^{1}}
\notag \\
\lesssim& (T^{\frac{1}{7}}M^{\frac{3}{5}}C_{0}+\eta) C_{0}\eta.\notag
\end{align}
Under the assumption that $\eta C_{0}\ll 1$, we choose
$T=T(C_{0},M,\eta)$ sufficiently small so that, by
\eqref{equ:low frequency,X1} and \eqref{equ:high frequency,X1},
\begin{align}
\bbn{\int_{0}^t e^{i(t-s)\Delta}F(u(s))ds
}_{X^{1}}\leq \eta,
\end{align}
which suffices to show that $\mathcal{T}_{f}$ maps $B(C_{0},\eta)$ into itself.

Next, we verify the contraction property.
For $u,v$ in this closed set $B(C_{0},\eta)$, applying \eqref{equ:trilinear,high,low} with $\sigma=1$, we have
\begin{align}\label{equ:contraction}
\bbn{\int_{0}^{t} e^{i(t-s)\Delta}
\lrs{F(u(s))-F(v(s))}ds}_{X^{1}}
\lesssim &\lrs{\eta C_{0}+T^{\frac{1}{7}}M^{\frac{3}{5}}C_{0}^{2}}
\n{u-v}_{X^{1}}\\
\leq &\frac{1}{2}\n{u-v}_{X^{1}}.\notag
\end{align}
Therefore, by the contraction mapping principle, we establish the local well-posedness in the space $B(C_{0},\eta)$.

We are left to prove \eqref{equ:stability}. Using \eqref{equ:trilinear,high,low} with $\sigma=-1$, we derive
\begin{align}\label{equ:difference,phi,f,g}
\n{\Phi(\cdot)f-\Phi(\cdot)g}_{X^{-1}}
\leq&\n{f-g}_{H^{-1}}+\frac{1}{2}\n{\Phi(\cdot)f-\Phi(\cdot)g}_{X^{-1}},
\end{align}
which implies \eqref{equ:stability}.

For the quintic case on $\T^{3}$, the corresponding high-low frequency estimate is given in \cite[Lemma 5.6]{CSZ22}, which reads
\begin{align*}
& \n{\int_{0}^{t} e^{i(t-s) \Delta} \prod_{j=1}^{5} \widetilde{v}_{j}(s) ds}_{X^\sigma} \\
 \lesssim&\n{v_{1}}_{X^\sigma} \prod_{j=3}^5\n{v_{j}}_{X^{1}}\lrs{T^{\frac{1}{12}}
 M^{\frac{1}{2}}\n{P_{\leq M} v_{2}}_{X^1}+\n{P_{>M} v_{2}}_{X^1}},
\end{align*}
where $\sigma\in\{-1,1\}$. Then, by repeating the proof for the cubic case, the quintic case follows analogously.
\end{proof}
\begin{remark}\label{remark:uniform stability}
Building on this local propagation estimate \eqref{equ:local propagation,strong flow} of the UTHT property, the uniform stability estimate \eqref{equ:stability} under the strong solution flow is then reduced to the high-low frequency multilinear estimates \eqref{equ:trilinear,high,low}.
The argument of Lemma \ref{lemma:stability,strong flow,UTFLS} extends to a broad range of cases, including the defocusing and focusing cases, as well as other powers and regularity indices whenever the
corresponding high-low nonlinear estimates, low-frequency bounds, and local strong solution
theory are available.
\end{remark}

\begin{lemma}\label{lemma:weak solution,UTFLS,strong flow}
Let $T_{0}>0$ and let $u$ be a weak solution to \eqref{equ:NLS} such that
\begin{align*}
u\in C([0,T_{0}];H_{\mathrm{weak}}^{1})
\end{align*}
and
\begin{align}\label{equ:weak,bounds,ut}
\sup_{t\in[0,T_{0}]}\n{u(t)}_{H^{1}}\leq C_{0},\quad
\sup_{t\in[0,T_{0}]}\n{P_{>M}u(t)}_{H^{1}}\leq\eta(C_{0}),
\end{align}
with the threshold $\eta(C_{0})\ll C_{0}^{-(p-2)}$. Then,
there exists $T=T(C_{0},M,\eta)>0$ such that, for every $t_{0}\in[0,T_{0}]$,
\begin{align}\label{equ:flow,weak solution,UTFLS}
\sup_{\tau\in[0,T]}\n{\Phi(\tau)u(t_{0})}_{H^{1}}\leq 2C_{0},  \quad
\sup_{\tau\in[0,T]}\n{P_{>M}\Phi(\tau)u(t_{0})}_{H^{1}}\leq 2\eta(C_{0}).
\end{align}

Moreover, for $t_{0}+t_{1}\in[0,T_{0}]$, $t_{1}\in[0,T]$, $t\in [0,T]$, it holds that
\begin{align}\label{equ:strong flow,weak solution,stabilty}
\n{\Phi(t)u(t_{0}+t_{1})-\Phi(t)\Phi(t_{1})u(t_{0})}_{H^{-1}}
\leq 2\n{u(t_{1}+t_{0})-\Phi(t_{1})u(t_{0})}_{H^{-1}}.
\end{align}
\end{lemma}
\begin{proof}
From \eqref{equ:weak,bounds,ut}, it follows that
$u(t_{0})\in \mathcal{S}(C_{0},M,\eta)$ for $t_{0}\in [0,T_{0}]$.
Using Lemma \ref{lemma:stability,strong flow,UTFLS}, we obtain $T_{1}(C_{0},M,\eta)$ such that
\begin{align}\label{equ:flow,weak solution,UTFLS,proof}
\sup_{\tau\in[0,T_{1}]}\|\Phi(\tau)u(t_{0})\|_{H^{1}}\leq2C_{0},
\
\sup_{\tau\in[0,T_{1}]}\|P_{>M}\Phi(\tau)u(t_{0})\|_{H^{1}}\leq2\eta.
\end{align}

Then, by \eqref{equ:weak,bounds,ut} and \eqref{equ:flow,weak solution,UTFLS,proof}, we get
\begin{align*}
\max\lr{
\n{u(t_{0}+t_{1})}_{H^{1}},
\n{\Phi(t_{1})u(t_{0})}_{H^{1}}}\leq& 2C_{0}\\
\max\lr{\|P_{>M}u(t_{0}+t_{1})\|_{H^{1}},
\n{P_{>M}\Phi(t_{1})u(t_{0})}_{H^{1}}}\leq& 2\eta,
\end{align*}
for $t_{0}+t_{1}\in[0,T_{0}]$ and $t_{1}\in[0,T_{1}]$.
Since we have $2\eta \ll (2C_{0})^{-(p-2)}$, we apply
Lemma \ref{lemma:stability,strong flow,UTFLS} again with parameters $2C_{0}$, $M$, $2\eta$. Consequently, there exists $T_{2}(2C_{0},M,2\eta)$ such that
\begin{align*}
\|\Phi(t)u(t_{0}+t_{1})-\Phi(t)\Phi(t_{1})u(t_{0})\|_{H^{-1}}
\leq2\n{u(t_{1}+t_{0})-\Phi(t_{1})u(t_{0})}_{H^{-1}},
\end{align*}
for $t\in [0,T_{2}]$.
Finally, by choosing
$$T(C_{0},M,\eta)=\min\lr{T_{1}(C_{0},M,\eta),T_{2}(2C_{0},M,2\eta)},$$
 we complete the proof of \eqref{equ:flow,weak solution,UTFLS} and \eqref{equ:strong flow,weak solution,stabilty}.
\end{proof}

\begin{proof}[\textbf{Proof of Proposition $\ref{lemma:weak uniqueness,threshold}$ on $\T^{d}$}]
Let $u$ be a weak solution to \eqref{equ:NLS} such that
\begin{align*}
u\in W^{1,\infty}(0,T_{0};H^{-1})\cap C([0,T_{0}];H^{1}_{\mathrm{weak}}),
\end{align*}
and
\begin{align}\label{equ:weak,bounds,ut,proof}
\sup_{t\in[0,T_{0}]}\n{u(t)}_{H^{1}}\leq C_{0},\ \sup_{t\in[0,T_{0}]}\n{P_{>M}u(t)}_{H^{1}}\leq \eta(C_{0}).
\end{align}
For almost every $t\in[0,T_{0})$, we have
\begin{align}\label{equ:limit,ut,H-1}
\lim_{\tau \to0+}\frac{u(t+\tau)-u(t)}{\tau}=i\Delta u(t)-iF(u(t)) \quad\text{in $H^{-1}$}.
\end{align}

By Lemma \ref{lemma:weak solution,UTFLS,strong flow}, there exists $T=T(C_{0},M,\eta)\in\lrs{0,\min\lr{1,T_{0}}}$ such that, for every $t\in[0,T_{0}]$ and $\tau\in [0,T]$, there exists a local strong solution $\Phi(\tau)u(t)$. Then, we also have
\begin{align}\label{equ:limit,Phit,H-1}
\lim_{\tau \to 0+}\frac{\Phi(\tau)u(t)-u(t)}{\tau}=i\Delta u(t)-iF(u(t))
\quad\text{in $H^{-1}$}.
\end{align}
Hence, for almost every $t\in[0,T_{0})$, we obtain
\begin{align}\label{equ:limit,ut,Phit,H-1}
\lim_{\tau\to 0+} \bbn{\frac{u(t+\tau)-\Phi(\tau)u(t)}{\tau}}_{H^{-1}}=0.
\end{align}

Next, we fix $t_{0}\in(0,T]$ and consider the flow $\Phi(t_{0}-t)u(t)$ for $t\in [0,t_{0}]$.
By the $H^{-1}$ uniform stability estimate \eqref{equ:strong flow,weak solution,stabilty} in Lemma \ref{lemma:weak solution,UTFLS,strong flow}, we deduce
\begin{align}\label{equ:H-1,stability,estimate,Phi,u}
\n{\Phi(t_{0}-t-\tau)u(t+\tau)-\Phi(t_{0}-t)u(t)}_{H^{-1}}
&\leq 2\n{u(t+\tau)-\Phi(\tau)u(t)}_{H^{-1}}.
\end{align}
Since $u\in W^{1,\infty}(0,T_{0};H^{-1})$ and the local strong solutions satisfy the uniform $H^{1}$ bound \eqref{equ:flow,weak solution,UTFLS}, we have
\begin{align*}
\n{u(t+\tau)-\Phi(\tau)u(t)}_{H^{-1}}\lesssim \tau,
\end{align*}
which implies the Lipschitz continuity of the map $t\mapsto \Phi(t_{0}-t)u(t)$.
Dividing \eqref{equ:H-1,stability,estimate,Phi,u} by $\tau$ and using
estimate \eqref{equ:limit,ut,Phit,H-1}, we obtain
\begin{align*}
\bbn{\frac{\Phi(t_{0}-t-\tau)u(t+\tau)-\Phi(t_{0}-t)u(t)}{\tau}}_{H^{-1}}\leq&
2\bbn{\frac{u(t+\tau)-\Phi(\tau)u(t)}{\tau}}_{H^{-1}}\to 0,
\end{align*}
which implies that
\begin{align*}
\frac{d}{dt}\bigl[\Phi(t_{0}-t)u(t)\bigr]=0\quad \text{in $H^{-1}$},
\end{align*}
for almost every $t\in[0,t_{0}]$. Consequently, we arrive at
\begin{align*}
u(t_{0})=\Phi(0)u(t_{0})=\Phi(t_{0})u(0), \quad 0\leq t_{0}\leq T.
\end{align*}

Since the local time $T(C_{0},M,\eta)$ is a fixed length, the same
argument can be applied on successive intervals of length $T$. Iterating
this procedure finitely many times, we extend the uniqueness to the
entire interval $[0,T_{0}]$. In particular, we obtain
\begin{align}\label{equ:identity,entire}
u(t)=\Phi(t)u(0),\quad t\in [0,T_{0}],
\end{align}
which suffices to conclude the uniqueness within the weak solution class under the condition \eqref{equ:weak,bounds,ut,proof}.
\end{proof}

\subsection{The Euclidean domains}
\begin{proof}[\textbf{Proof of Proposition $\ref{lemma:weak uniqueness,threshold}$ on $\R^{d}$}]
We suppose that $u$ and $v$ are two weak solutions on $[0,T_{0}]$, with $u(0)=v(0)$, and assume that
\begin{equation}\label{equ:u,v,utht,proof}
\left\{
\begin{aligned}
\max \lr{\n{u}_{L^{\infty}([0,T_{0}]; H^{1})},
\n{v}_{L^{\infty}([0,T_{0}];H^{1})}} & \leq C_{0}, \\
\max \lr{\n{P_{>M} u}_{L^{\infty}([0,T_{0}]; H^{1})},\n{P_{>M}v}_{L^{\infty}([0,T_{0}]; H^{1})}}& \leq \eta .
\end{aligned}
\right.
\end{equation}
We take a difference $w=u-v$ and define the Strichartz norm
$$
\n{w}_{\mathcal{S}[0,T]}=\n{w}_{L^{\infty}\left([0,T] ; L^2\right)}+\n{w}_{L^2\left([0,T] ; L^r\right)} .
$$
where
\begin{align*}
r=\frac{2d}{d-2}, \quad r^{\prime}=\frac{2d}{d+2}.
\end{align*}
We notice that
\begin{align*}
w(t)=&\pm i\int_{0}^{t}e^{i(t-s)\Delta}(|v+w|^{p-1}(v+w)(s)-|v|^{p-1}v(s))ds\\
=&\int_{0}^{t}e^{i(t-s)\Delta}\lrc{\lrs{A_{high}+A_{low}}w+\lrs{B_{high}+B_{low}}\ol{w}}ds,
\end{align*}
where $A_{\text{high}}$ and $B_{\text{high}}$ collect the coefficients carrying at least one high-frequency factor, and $A_{\text{low}}$ and $B_{\text{low}}$ consist of all low-frequency coefficients.

Based on estimates \eqref{equ:u,v,utht,proof},
it is standard to check that
$$
\n{A_{\text{high}}}_{L^{\infty}([0,T];L^{\frac{d}{2}})} +
\n{B_{\text{high}}}_{L^{\infty}([0,T];L^{\frac{d}{2}})}\lesssim \eta C_{0}^{p-2},
$$
and
$$
\n{A_{\text{low}}}_{L^{\infty}([0,T];L^{\infty})}+
\n{B_{\text{low}}}_{L^{\infty}([0,T];L^{\infty})} \lesssim M^{2} C_{0}^{p-1}.
$$
Applying the retarded endpoint Strichartz estimate \cite[Theorem 1.2]{KT98} to the high-frequency terms in $L^{2}([0,T];L^{r^{\prime}})$and to the low-frequency terms in
$L^{1}([0,T];L^{2})$, we arrive at
\begin{align*}
&\n{w}_{\mathcal{S}([0,T])} \\
\lesssim& \n{A_{high}w}_{L^{2}([0,T];L^{r^{\prime}})}+
\n{B_{high}w}_{L^{2}([0,T];L^{r^{\prime}})}+\n{A_{low}w}_{L^{1}([0,T];L^{2})}
+\n{B_{low}w}_{L^{1}([0,T];L^{2})}\\
\lesssim&(\eta C_{0}^{p-2}+T M^{2} C_{0}^{p-1})\n{w}_{\mathcal{S}([0,T])}.
\end{align*}

Choosing $\eta C_{0}^{p-2}$ sufficiently small and then $T$ sufficiently small, we deduce that $w=0$. A standard iteration argument then yields uniqueness in the UTHT class.
As shown in \cite{CH19,CH22unc}, a strong solution lies in the UTHT class.
Combining this with strong-weak uniqueness and standard continuation, we obtain the representation \eqref{equ:identity,entire,prop} of the strong solution flow.

\end{proof}
\section{Global weak uniqueness}\label{section:Global weak uniqueness}
In this section, we study the defocusing energy-critical NLS and reformulate the UTHT criterion as an energy defect threshold criterion. Along the strong solution flow, energy conservation propagates the EDT condition for the weak solution and, combined with the global well-posedness theory for the strong solution, ultimately gives global weak uniqueness.
\subsection{Energy defect threshold criterion}\label{section:Local regularity below the energy defect threshold}
For $t_{0}\in I$, we define the energy defect
\begin{align}\label{equ:defect}
D(u,t_{0})=\lim_{\ve \to 0}\operatorname*{ess\, sup}_{|t-t_{0}|<\ve}E(u(t))-E(u(t_{0})).
\end{align}
Since $E(t)$ is lower semicontinuous, its supremum and essential supremum coincide on every open subinterval.
We note that, by the $H^{1}$ weak continuity of $u(t)$, the quantity $D(u,t_{0})$ is well-defined and nonnegative.

Next, we present a uniqueness criterion for the class of weak solutions satisfying the EDT condition.
\begin{lemma}[EDT criterion]\label{lemma:energy defect,strong-weak}
Let $I_{0}$ be an open interval, $t_{0}\in I_{0}$, and $u\in C(I_{0};H_{\mathrm{weak}}^{1})$ be a weak solution to \eqref{equ:NLS} with an $H^{1}$ energy bound $C_{0}\geq 1$, and set
$\delta(C_{0})=\frac{\eta(C_{0})^{2}}{16}$, where $\eta(C_{0})$ is the threshold from Proposition
$\ref{lemma:weak uniqueness,threshold}$. If $D(u,t_{0})<\delta(C_{0})$, then there exists
$\varepsilon>0$ such that
\begin{align}\label{equ:ut,phit,t0}
u(t)=\Phi(t-t_{0})u(t_{0}),\quad \f t\in (t_{0}-\ve,t_{0}+\ve)\subset I_{0}.
\end{align}
\end{lemma}
\begin{proof}
We begin with the identity
\begin{align*}
&\frac{1}{2}\n{\nabla u(t)-\nabla u(t_{0})}_{L^{2}}^{2}\\
=&E(u(t))-E(u(t_{0}))-\operatorname{Re}\lra{\nabla u(t)-\nabla u(t_{0}),\nabla u(t_{0})}+\frac{1}{p+1}
\lrs{\n{u(t_{0})}_{L^{p+1}}^{p+1}-\n{u(t)}_{L^{p+1}}^{p+1}}.\notag
\end{align*}
By the weak $H^{1}$ continuity of $u(t)$, we have
\begin{align*}
&\lim_{t\to t_{0}}\lra{\nabla u(t)-\nabla u(t_{0}),\nabla u(t_{0})}=0,\\
&\n{u(t_{0})}_{L^{p+1}}\leq \liminf_{t\to t_{0}}\n{u(t)}_{L^{p+1}}.
\end{align*}
Consequently, together with the $L^{2}$ continuity obtained by interpolation, we obtain
\begin{align*}
\limsup_{t\to t_{0}}\n{u(t)-u(t_{0})}_{H^{1}}^{2}\leq 2D(u,t_{0}),
\end{align*}
which, combined with the condition that $D(u,t_{0})<\delta=\frac{\eta^{2}}{16}$, implies that there exists $\ve>0$ such that
\begin{align}\label{equ:ut,ut0,D}
\sup_{|t-t_{0}|<\ve}\n{u(t)- u(t_{0})}_{H^{1}}^{2}\leq \frac{\eta^{2}}{4}.
\end{align}
For this threshold $\eta$ and fixed time $t_{0}$, we choose $M$ such that
\begin{align}\label{equ:PM,ut0,eta}
\n{P_{>M}u(t_{0})}_{H^{1}}\leq \frac{\eta}{2}.
\end{align}

Using \eqref{equ:ut,ut0,D} and \eqref{equ:PM,ut0,eta}, we deduce
\begin{align}\label{equ:PM,ut,threshold}
\sup_{|t-t_{0}|<\ve }\n{P_{>M}u(t)}_{H^{1}}\leq \n{P_{>M}u(t_{0})}_{H^{1}}+\sup_{|t-t_{0}|<\ve }\n{ u(t)-u(t_{0})}_{H^{1}}\leq \eta.
\end{align}

Applying \eqref{equ:identity,entire,prop} in Proposition
\ref{lemma:weak uniqueness,threshold}, after time translation and in both time directions, we obtain \eqref{equ:ut,phit,t0}.
\end{proof}

\subsection{Strong-weak uniqueness}\label{section:Strong-weak uniqueness}
\begin{proof}[Proof of Theorem $\ref{theorem:weak uniqueness}$]
Let $I_{0}\Subset I$ be a bounded open interval. By the local-in-time $H^{1}$ energy bound, there exists $C_{0}\geq 1$ such that
\begin{align}
\n{u}_{L_{t}^{\infty}(I_{0};H^{1})}\leq C_{0}<\infty,
\end{align}
and hence
\begin{align}
\operatorname*{ess\, sup}_{t\in I_{0}}E(u(t))<\infty.
\end{align}

Let $\delta=\delta(C_{0})$ be the energy defect threshold given by Lemma \ref{lemma:energy defect,strong-weak}. Since the energy can be made arbitrarily close to its essential supremum, we may choose a time $t_{0}\in I_{0}$ such that
\begin{align*}
E(u(t_{0}))>\operatorname*{esssup}_{t\in I_{0}}E(u(t))-\delta(C_{0}),
\end{align*}
which implies the energy defect threshold property at time $t_{0}$, that is,
\begin{align}\label{equ:Dut0,delta}
D(u,t_{0})\leq \operatorname*{esssup}_{t\in I_{0}}E(u(t))-E(u(t_{0}))<\delta.
\end{align}
On the other hand, the well-known global theory \cite{Bou99,CKSTT08,HTT11,IP12,RV07,Kwak26,Yue21} provides a global strong solution
\begin{align*}
v(t)=\Phi(t-t_{0})u(t_{0}),\quad t\in \R.
\end{align*}

By Lemma \ref{lemma:energy defect,strong-weak}, there exists $\ve_{0}>0$ such that
\begin{align}\label{equ:ut=vt,t0,ve0}
u(t)=v(t),\quad \f t\in(t_{0}-\ve_{0},t_{0}+\ve_{0}),
\end{align}
which implies local strong-weak uniqueness near $t_{0}$.
We notice that the strong solution $v(t)$ conserves energy, namely,
\begin{align}\label{equ:energy conservation,vt,ut0}
E(v(t))=E(u(t_{0})),\quad t\in \R.
\end{align}
Therefore, $u(t)$ in a neighborhood of $t_{0}$ still satisfies the energy defect threshold condition.
Based on this observation, we can employ a standard bootstrap argument to extend the local uniqueness to global uniqueness.

To proceed, we set
\begin{align*}
S=\lr{t\in I_{0}:u(t)=v(t)}.
\end{align*}
From \eqref{equ:ut=vt,t0,ve0}, the set $S$ is nonempty, and it is closed by the weak $H^{1}$ continuity of $u(t)$ and $v(t)$. It remains only to prove that $S$ is open.
For every $t_{1}\in S$, combining $u(t_{1})=v(t_{1})$, the energy conservation \eqref{equ:energy conservation,vt,ut0}, and the small energy defect estimate \eqref{equ:Dut0,delta}, we get
\begin{align*}
D(u,t_{1})\leq& \operatorname*{esssup}_{t\in I_{0}}E(u(t))-E(u(t_{1}))\\
=&\operatorname*{esssup}_{t\in I_{0}}E(u(t))-E(v(t_{1}))\\
=&\operatorname*{esssup}_{t\in I_{0}}E(u(t))-E(u(t_{0}))<\delta
\end{align*}
Applying Lemma \ref{lemma:energy defect,strong-weak} again, we obtain
\begin{align*}
u(t)=v(t),\quad \f t\in(t_{1}-\ve_{1},t_{1}+\ve_{1}).
\end{align*}
This shows that $S$ is also open. Therefore, we conclude that $S=I_{0}$.

Finally, combining the strong-weak uniqueness established for bounded open intervals with the global existence of the strong solution, the uniqueness can be extended to the entire interval $I$. In other words, within the class of weak solutions, we obtain
\begin{align*}
u(t)=\Phi(t-t_{0})u(t_{0}),\quad \f t\in I.
\end{align*}
 Hence, we complete the proof of Theorem \ref{theorem:weak uniqueness}.

\end{proof}

\section{Uniqueness for the energy-critical GP hierarchy}\label{section:Uniqueness for the Energy-critical GP hierarchy}

In this section, we apply Theorem \ref{theorem:weak uniqueness} to prove uniqueness for the
energy-critical GP hierarchy.
We use the uniqueness principle of Ammari, Liard, and Rouffort \cite[Theorem 1.4 and its proof]{ALR20} in the following form adapted to the present setting. \footnote{Although \cite[Theorem 1.4]{ALR20} is formulated for open intervals, it remains valid on $[0,T]$ under the continuity assumptions stated here. In fact, two hierarchies can be extended identically to $(-\ve,0]$ using the strong NLS flow, so the \cite[Theorem 1.4]{ALR20} for open intervals applies with $t=0$ as an interior initial time. For simplicity, we omit the proof.
 }

\noindent \textbf{Uniqueness principle.}
Suppose that $v:H^{1}(\Lambda)\to H^{-1}(\Lambda)$ satisfies (A0) below.
Then uniqueness of weak solutions on $[0,T]$ to
\begin{align*}
\partial_{t}u=v(u),\quad u(0)=u_{0},
\end{align*}
in $L^{\infty}(0,T;H^{1})\cap W^{1,\infty}(0,T;H^{-1})$
implies uniqueness of solutions to the associated symmetric hierarchy
satisfying (A1), for each initial hierarchy.
Here, the assumptions are as follows.

\begin{enumerate}[wide=0pt, label=(\arabic*), align=left]
\item[(A0)]
The vector field $v$ is Borel and bounded on bounded sets, and
$v$ is $U(1)$-invariant, namely,
\begin{align*}
v(e^{i\theta}u)=e^{i\theta}v(u),\quad \theta\in\R.
\end{align*}

\item[(A1)]
The hierarchy satisfies the following three conditions.
\begin{enumerate}[wide=0pt, label=(\arabic*), align=left]
\item[(i)]
 For every $t\in [0,T]$, there is a $U(1)$-invariant Borel probability
measure $\mu_{t}$ on the closed unit ball of $L^{2}(\Lambda)$ such that
\begin{align*}
\gamma^{(k)}(t)
=\int (|\phi\rangle\langle\phi|)^{\otimes k}\,d\mu_{t}(\phi),
\quad k\geq 1.
\end{align*}

\item[(ii)]
There is a constant $C_{0}>0$, independent of $k$, such that
\begin{align*}
\sup_{t\in [0,T]}
\|S^{(1,k)}\gamma^{(k)}(t)S^{(1,k)}\|_{\mathcal L_{k}^{1}}
\leq C_{0}^{2k},\quad k\geq 1.
\end{align*}

\item[(iii)]
For every $k\geq 1$, the map
\begin{align*}
t\mapsto S^{(-1,k)}\gamma^{(k)}(t)S^{(-1,k)}
\end{align*}
is weak-$*$ continuous in $\mathcal L_{k}^{1}$ on $[0,T]$, where the weak-$*$
topology is induced by testing against compact operators.
\end{enumerate}
\end{enumerate}

\begin{proof}[Proof of Theorem~\ref{thm:hierarchy}]
With $b_{0}=1$, we take
\begin{align*}
v(u)=i\Delta u-i|u|^{4}u.
\end{align*}
By the Sobolev embedding $H^{1}(\Lambda)\hookrightarrow L^{6}(\Lambda)$,
the vector field $v$ is continuous from $H^{1}$ to $H^{-1}$ and satisfies
\begin{align*}
\n{v(u)}_{H^{-1}}
\lesssim \n{u}_{H^{1}}+\n{u}_{H^{1}}^{5}.
\end{align*}
Together with $U(1)$-invariance, this verifies (A0).
Conditions (1), (2), and (4) of Definition \ref{def:weak solution,gp} yield (A1),
and condition (3) is the integral formulation of the symmetric
hierarchy associated with $v$.

Theorem \ref{theorem:weak uniqueness} provides the required uniqueness of weak NLS solutions.
Applying the above uniqueness principle, we conclude that weak solutions
of \eqref{equ:gp hierarchy,quintic} with the same initial hierarchy
coincide on $[0,T]$.
\end{proof}

\noindent\textbf{Disclosure of AI assistance.}
The key idea of the weak uniqueness within the UTHT class originates from the previous work \cite{CH19,CH22unc}. ChatGPT-6 Astra (OpenAI) was used to assist in exploring potential simplifications of the UTHT criterion. Based on our understanding of the UTHT criterion, we work out this simplified proof of Proposition \ref{lemma:weak uniqueness,threshold}.
 The authors take full responsibility for the accuracy, originality, and entire content of the manuscript.

\noindent\textbf{Data Availability Statement.}
Data sharing is not applicable to this article as no datasets were generated or analysed during the current study.

\noindent\textbf{Conflict of Interest.}
The authors declare that they have no conflict of interest.

\bibliographystyle{abbrv}
\bibliography{references}

\end{document}